# Application of substring searching methods
# to group presentations[*]

George Havas[†] and Mark Ollila
Key Centre for Software Technology
Department of Computer Science
University of Queensland
Queensland 4072
Australia

## Abstract

An important way for describing groups is by finite presentations. Large presentations arise in practice which are poorly suited for either human or computer use. Presentation simplification processes which take "bad" presentations and produce "good" presentations have been developed. Substantial use is made of substring searching and appropriate techniques for this context are described. Effective use is made of signatures and change flags. Change flags are shown to be the most beneficial of the methods tested here, with very significant performance improvement. Experimental performance figures are given.

## 1 Introduction

Finitely presented groups have been much studied. All requisite mathematical background is provided in [14, Chapter 1]. A comprehensive book on computation with finitely presented groups by Charles Sims [20] is soon to appear. An overview of algorithms for such groups is included in [3].

A finitely presented group may be given by a presentation

$$G = \langle g_1, \ldots, g_d \mid R_1, \ldots, R_n \rangle$$

where the $g_i$ are generators and the $R_j$ are relators. Generally speaking, presentations are good if they are short: few generators; and few relators of reasonable length. This makes them relatively intelligible to humans and also often makes them better suited for computer calculations. We are interested in the situation where we have what we regard as a bad presentation for a group and we wish to find a good presentation. This kind of situation may arise combinatorially in a number of ways. The most common is through subgroup presentation procedures.

Based on theorems of Reidemeister and Schreier, various algorithms for computing presentations of subgroups of finite index in finitely presented groups have been developed. These include Havas [7], Arrell and Robertson [2] and a reduced Reidemeister-Schreier program available in SPAS [4] which is based on the description by Neubüser [17]. Subgroup presentations produced via Reidemeister-Schreier processes generally have large numbers of generators and relators and are not well-suited for either human or computer use.

A theorem of Tietze proves that, given two presentations of a group, a sequence of transformations exists which demonstrates that the presentations are of the same group. However there is no general algorithm for finding such a sequence, a consequence of unsolvability results in group theory. Tietze transformation procedures which input a "bad" presentation and output a "good" presentation for a group have been developed (Havas, Kenne, Richardson and Robertson [10], Robertson [18]). New Tietze transformation procedures written in the higher level GAP [19] language (with special kernel support) are under development by Volkmar Felsch and Martin Schönert at Aachen. We study in particular the substring searching part of such Tietze transformation processes.

## 2 Presentation simplification

Given a bad presentation for a group there are a number of different approaches to presentation simplification aimed at producing a better presentation. Charles Sims [20, Chapter 6.4] discusses

---

use of Knuth-Bendix and coset enumeration based methods in addition to the Tietze transformation methods which we study. The effectiveness of Tietze transformation programs is witnessed in applications such as Newman and O'Brien [16] and Havas and Robertson [9], where some finitely presented groups are proved to be soluble and their maximal soluble quotients are computed.

Any group relator is equivalent to the relator obtained from it by free reduction. It is also equivalent to all cyclic permutations of itself and its formal inverse. This leads to a natural canonical form for relators, namely the least of all of these equivalents (in length by lexicographic order). Efficient methods for finding such canonical forms are discussed by Iliopoulos and Smyth [12]. Using this canonical form, group presentations may be stored as a sorted sequence of canonical representatives of relators, which facilitates the discovery (and removal) of duplicates.

Three main principles used by Tietze transformation methods to simplify presentations [10] are: short eliminations; substring replacements; long eliminations.

- *Short eliminations.*
  All relators of length 1 and non-involutory relators of length 2 are used to eliminate generators and their associated relators.

- *Substring replacements.*
  Relators are shortened by replacing long substrings by shorter equivalent strings. First substring searching is performed. A relator $R_i$ is chosen and other relators are searched for a matching substring $v$ in a cyclic permutation $uv$ of $R_i$ or its inverse and in a cyclic permutation $wv$ of $R_j$, with the length of $v$ greater than the length of $u$. Then, when such a match is found, the relator $R_j$ is replaced by a canonical representative of the shorter, equivalent relator $wu^{-1}$. One substring replacement pass involves the application of this process with $R_i$ running once through all relators in the presentation and being used to search all subsequent relators.

- *Long eliminations.*
  Redundant generators (generators which occur only once in some relator) and their associated relators are eliminated using relators with length greater than 2.

Tietze transformation programs function by working through these steps in some sensible order, guided by heuristics. In automatic mode, short eliminations are performed till no more are possible. Then substring replacement passes are made. If possible, short eliminations are again done, otherwise long eliminations are performed. This sequence is repeated until no further improvements to the presentation are possible this way. Short eliminations and substring replacements reduce the total length of the presentation. Long eliminations can, and generally do, increase the length, often quite significantly. The substring searching component of substring replacements is by far the most time consuming part of presentation simplification, which is why we focus on it here.

## 3 An easy example

To clarify the presentation simplification process we give an easy example commencing with 9 generators and 9 relators. Our aim is to minimize the number of generators and relators and the length of the relators. Consider the Fibonacci group $F(2,9) =$
$$\langle x_1, \ldots, x_9 \mid x_1 x_2 x_3^{-1}, \ldots, x_8 x_9 x_1^{-1}, x_9 x_1 x_2^{-1} \rangle,$$
which provides one of our main examples later. For convenience we represent $x_1$ by $a$ and $x_1^{-1}$ by $A$, $x_2$ by $b$, ..., $x_9$ by $i$ and $x_9^{-1}$ by $I$. For such a small set of relators it is more convenient neither to use canonical forms nor to sort the relators into length by lexicographical order, since this makes it easier to follow the working. Here we underline useful common substrings in sets of relators to make them readily visible.

Written this way our initial set of relators is $\{abC, bcD, cdE, deF, efG, fgH, ghI, hiA, iaB\}$. If we choose to eliminate $c = ab$, $f = de$ and $i = gh$ in one batch the next set of relators is $\{babD, abdE, edeG, degH, hghA, ghaB\}$. With no useful common substrings in this set, we choose to eliminate $d = bab$ next. This gives $\{abbabE, \underline{ebabe}G, \underline{babeg}H, hghA, ghaB\}$.

Now $R_2$ and $R_3$ have a useful common substring, $babe$, which is replaced in $R_3$ by its equivalent value $Eg$ from $R_2$. Thus the set of relators becomes $\{abbabE, ebabeG, EggH, hghA, ghaB\}$. This set has no additional useful common substrings, so we continue by eliminating $h = GbA$, giving $\{abbabE, ebabeG, EggaBg, GbAbAA\}$, with no useful substrings. Eliminating $g = bAbAA$ we get $\{abbabE, \underline{ebabe}aaBaB, \underline{EbAbAA}bAbAAbAA\}$. In this case a permutation of the inverse of $R_2$ and $R_3$ now have a useful common substring $EbAbAA$, so we replace it in $R_3$ by its equivalent from $R_2$, $babe$, giving a new set $\{abbabE, ebabeaaBaB, babebAbAAbAA\}$. There are no other useful common substrings, and we continue by eliminating $e = abbab$, giving $\{abbabbababbabaaBaB, bababbabbAbAAbAA\}$, which is not further simplified and is, with relators in the other order, the presentation we use in the next section.

In this simple example substring replacements had some beneficial effect, with the final relators being a little shorter than they otherwise would have been. This is greatly magnified in larger examples. To make matters worse in large examples, the number of eliminations possible may be substantially reduced if substring replacements are not made.

## 4   Substring searching

Aho [1] and Gonnet and Baeza-Yates [5, Chapter 7] present algorithms and data structures for substring searching in various situations. However the case considered here differs substantially from those covered there. The major distinguishing features of our situation are: all strings are (in effect) circular; formal inverses are (implicitly) present; many substrings are simultaneously sought; the text is dynamic, changing very often. We study effective methods for this presentation simplification environment.

In spite of the theoretical worst case inferiority of brute force searching, its average case performance is linear in the length of the text being searched. Furthermore, Gonnet and Baeza-Yates [5, Table 7.4] show that it performs quite well in practice. We use a variant of brute force searching which enables us to search for many strings simultaneously at no extra cost. Improvements available through faster string searching algorithms like those of Hume and Sunday [11] will be the subject of future work, as will utilization of Karp and Rabin [13] type methods. No effort has yet been made to investigate substring searching with finite automata in this context, however the dynamic nature of the relators being searched is likely to be a significant impediment.

We first look at how substring searching is done for one pair of relators, $R_1 < R_2$ (in the length by lexicographic order). In order to shorten $R_2$ any useful common substring must have length greater than half the length of $R_1$. This means that it will contain either the first symbol of $R_1$ or a middle symbol, or the inverse of one of those. (Further, if $R_1$ is a nontrivial power, a useful substring must contain the first symbol or its inverse.) The inverse of a symbol is represented by the negative of the symbol's representation. Thus we commence by searching for one of at most two absolute values as starting points in $R_2$. When a match of absolute values is found we try to extend the match circularly both backwards and forwards until it is long enough to be useful.

A substring replacement pass potentially has beneficial effect after any alteration to the presentation. Thus substring replacement passes are made after generator eliminations and after successful substring replacement passes, that is, passes which result in a change to any relator. Various different strategies are described in [10] for deciding when to perform substring replacement passes, in order to reduce time spent in relatively useless substring searches.

We exemplify the performance of three improved methods for substring searching applied to group presentations by considering two specific examples in detail. The performance gains demonstrated here typify the improvements achieved in this application area by these methods.

We study two subgroup presentation applications, looking at performance on presentations $\mathcal{J}$ and $\mathcal{F}$. Presentation $\mathcal{J}$ is of the index 100 subgroup $\langle a, b, b^{ca^{-1}c} \rangle$ in the Janko simple group

$J_2 = \langle a, b, c \mid a^3, b^3, c^3, aba(bab)^{-1}, (ca)^5, (cb)^5,$
$(cb^{-1}cb)^2, a^{-1}baca^{-1}ba(bac)^{-1}ac^{-1},$
$aba^{-1}cab(ca)^{-1}ab^{-1}(ca)^{-1} \rangle.$

Presentation $\mathcal{J}$ was obtained using the standard Reidemeister-Schreier method described in [7] and has 201 generators, 510 relators with longest relator of length 12, and total relator length 2795. Presentation $\mathcal{F}$ is of the index 152 subgroup $\langle [a^2, b] \rangle^G$ in the Fibonacci group

$G = F(2, 9) = \langle a, b \mid$
$babab^2ab^2a^{-1}ba^{-2}ba^{-2}, ab^2ab^2abab^2aba^2b^{-1}ab^{-1} \rangle.$

$\mathcal{F}$ also was obtained using the program described in [7] and plays a crucial role in proving $F(2, 9)$ to be infinite, see [8, 15]. Presentation $\mathcal{F}$ has 153 generators, 304 relators with longest relator of length 13, and total relator length 2119.

Using specific default strategies for simplification, presentation $\mathcal{J}$ is reduced by Tietze transformations to a presentation with 3 generators, 43 relators with longest relator of length 21, and total relator length 504. This presentation can be comfortably used for other computations such as coset enumeration, for example. Likewise, presentation $\mathcal{F}$ is reduced to a presentation with 21 generators, 108 relators with longest relator of length 301, and total relator length 14510. Even though the total length of this presentation is greater than that of the initial presentation, it is more useful in the sense that it directly implies infiniteness (see [15]), a result by no means obvious from the initial presentation.

Our performance figures all refer to exactly the same sequences of elimination and replacement steps, with only the substring searching part changed. In both examples, substring replacement passes were made periodically after changes to the presentation. In particular, the passes were made after batches of short eliminations, after individual long eliminations, and after any substring replacement pass which led to change. This is a relatively expensive strategy which leads to more substring searching than other more conservative strategies. It is thus well suited for testing substring searching methods.

The key statistics to note are the number of pairs of relators searched with the straightforward method and the success rate. In the case of $\mathcal{J}$, a total of 6413797 searches were made, with only 1907 successful. The total time taken (on an otherwise lightly loaded Sun SS2) was 156 seconds. For $\mathcal{F}$, a total of 8614464 searches were made, with only 2836 successful. The total time taken was 446 seconds.

# 5  Signatures

It is clear from the very low success rate of searches in these applications that we can achieve efficiencies if we can speed up or even avoid unsuccessful searches. In this section we show how most of the unsuccessful searches can be made substantially faster. The key to the time saving is the concept of signatures, introduced by Harrison [6]. The speed ups come from the replacement of some time consuming substring searches by much faster tests which reveal that no useful match is possible.

Harrison points out that a string can be characterized by its set of substrings. Further, the set of substrings of a particular length is a string invariant, in the sense that if $W_1$ is a substring of $W_2$ then the set of substrings of a given length of $W_1$ is a subset of the set of substrings of the same length of $W_2$. Direct computation with sets of substrings would be a complicated process. Harrison goes on to associate with each set of substrings a bit string computed by hashing. For sets of substrings of length $k$ he calls this the *hashed k-signature*. There are two important parameters associated with this: $k$; and $m$, the length of the signature. Harrison goes on to describe how $k$ and $m$ associated with the method may be sensibly established and analyzed. We choose substrings of length 2 ($k = 2$) and signatures of length 32 ($m = 32$) for reasons which become clear later.

Thus, the idea for checking whether one string is a substring of another is to first compute their signatures. Signatures $S_1$ and $S_2$ from words $W_1$ and $W_2$, where $W_1$ is not longer than $W_2$, are compatible if $S_1 \& S_2 = S_1$, where & represents the bitwise "and" operation. If the signatures are incompatible then we know that $W_1$ is not a substring of $W_2$. If the signatures are compatible this does not guarantee that there is a match. We must continue by performing a normal substring search to check whether or not there is an actual match. Provided the signature generation and testing are fast and that enough normal substring tests (which are slow) are avoided, we save time.

Now let us look at our situation. We check relators for matches in pairs, $R_1$ and $R_2$, with $R_1 < R_2$. If the length of $R_1$ is $l$ then there are $2l$ substrings of $R_1$ any of which could be minimal length useful matches, possibly all distinct. This is because a minimal length useful match has length $\lceil (l+1)/2 \rceil$, may commence at any symbol in the relator, and comes from either (a cyclic permutation of) the relator itself or its formal inverse. These are the patterns we need to look for.

Given that we wish to search all cyclic permutations of $R_2$ it is most desirable that our target signatures be invariant under cyclic permutation. This avoids the need to have different signatures for each permutation, which would multiply the amount of work required by the length of $R_2$. Consider $R = \prod_{i=1}^{n} x_i$. The standard set of length 2 substrings is $\{x_1 x_2, x_2 x_3, \ldots, x_{n-1} x_n\}$. To make it permutation invariant we simply add $x_n x_1$ to the set.

Even with permutation invariance, a straight-forward implementation of signatures would still involve $2l$ signatures for $R_1$ and one signature for $R_2$. Since inversion is needed for only one of $R_1$ and $R_2$ and we have to consider $l$ substrings of $R_1$ and only $R_2$ itself, it is better to apply the inversion to $R_2$. Instead of $2l$ signatures for the substrings of $R_1$ and its inverse and one signature for $R_2$, we can use $l$ signatures for $R_1$ and 2 signatures for $R_2$. If any of the signatures of $R_1$ turns out to be compatible with a signature for $R_2$ then we must do a normal substring search. However we do better than that. We also make our signatures invariant under formal inversion.

One way to make them inversion invariant is to add all of the formal inverses of the members to the set of substrings. This has the disadvantage of doubling the size of the set, which would reduce the distinguishing power when combined with "random" hashing to produce the signature. However, for length 2 substrings, we have a much better way. We choose an inversion invariant hashing function $h$. We define

$$h(x_1 x_2) = (|f(x_1) - f(x_2)|) \bmod m$$

where $f(x_i)$ is an injection into the positive integers and $f(x_i^{-1}) = -f(x_i)$. It follows immediately that $f(x_2^{-1} x_1^{-1}) = f(x_1 x_2)$, that is, $h$ is inversion invariant. This means that we halve the number of signature comparisons we need, from $2l$ to $l$.

The selection of $f$ is easy, based on the unique numerical value of the internal representation of $x_i$. The simplest choice for $m$ in the context of current machine architectures is 32. This allows us to store our hashed 2-signature in one machine word and to test for compatibility of signatures by simple logical operations readily available on most machines. Two signatures $S_1$ and $S_2$ (from relators $R_1$ and $R_2$ with $R_1 < R_2$) are compatible if and only if $S_1 \& S_2 = S_1$.

We implemented substring searching based on these methods and obtained the following results. For $\mathcal{J}$, 5193796 out of 6413797 pairs of relators previously tested were incompatible. This left only 1220001 pairs to test, a saving of 81%. The total time taken was 85 seconds, a saving of 46%. For $\mathcal{F}$, 7057646 out of 8614464 pairs of relators previously tested were incompatible. This left only 1556818 pairs to test, a saving of 82%. The time taken was 295 seconds, a saving of 34%. (In these timings we did not take advantage of the fact that specific signature compatibility implies a specific potential match to check. Doing so might improve the time saving gained from signatures. However we simply used complete normal substring testing for the pair of relators, as described in §4.)

Signature computation and comparison is by no means free. The time taken for this was 44 seconds for $\mathcal{J}$ and 139 seconds for $\mathcal{F}$. Thus, our saving of 71 seconds for $\mathcal{J}$ was made up of 115 saved in normal substring searches avoided less the extra 44 seconds on signatures. Likewise, for $\mathcal{F}$, there were 290 seconds of normal substring searches saved. However the space required by the signatures is negligible in our implementation. For each target, one word is reserved, the signature is computed initially, and it is recomputed only when the target changes. On the other hand, the $l$ signatures for the patterns are recomputed as required, to avoid doubling the amount of space required for the presentation. This is the usual time/space trade-off. Saving all of these signatures would avoid much recomputation, but at high storage cost.

Further analysis revealed that with signature length $m = 32$ the compatibility test lost its effectiveness as the Tietze transformation process progressed. Two things caused this. The number of different generators was reducing and the average length of the relators was increasing. Both of these reduce the chance of a pair of relators being incompatible via a signature test. It is possible that this could be addressed by increasing $m$, however this does not appear worthwhile in view of the results with change flags.

# 6 Change flags

Signatures provide a way of speeding up most unsuccessful searches. However we do much better by avoiding unnecessary searches.

The sequence of relators which make up the presentation changes very frequently. Eliminations (short and long) and succcessful replacement passes make changes to relators. However, by no means all relators are changed between substring replacement passes. It is futile to embark on substring searches with a pair of relators neither of which has changed since the previous replacement pass. No useful common substrings were found on the last pass (since the lexicographically later one is unchanged), so none will be found now.

However the naive description of §2 and also early implementations of Tietze transformation programs blithely compare each relator $R_i$ with every subsequent relator in the relator sequence. The simple way to avoid the futile searches is to flag relators which have changed since the last substring replacement pass (when all relators pairs were searched). A pair of relators is then searched only if at least one of the pair is flagged. Thus, two flags are associated with each relator: *cslp* (changed since last pass); *ctp* (changed this pass). Initially the *cslp* flag is set for all relators and the *ctp* flag unset. During both short and long eliminations *cslp* flags are set for all relators which are changed. During substring replacement passes the *ctp* flag is set for any relator which changes. At the end of a substring replacement pass the *cslp* flags are (logically, not physically) reset to the current *ctp* flags and all *ctp* flags are reset.

For $\mathcal{J}$, 5960294 out of 6413797 pairs of relators previously tested were unchanged. This left only 453503 pairs to test, a saving of 93% over normal testing and a reduction of 63% on signature methods. The time taken was 32 seconds, a saving of 79% on normal methods and 62% as against signatures. For $\mathcal{F}$, 7843972 out of 8614464 pairs of relators previously tested were unchanged. This left only 770492 pairs to test, a saving of 91% over normal testing and a reduction of 51% on signature methods. The time taken was 108 seconds, a saving of 76% on normal methods and 63% as against signatures.

Notice that there is substantially less computation involved in this flag setting and testing than the work done in signature based tests. The time taken for flag processing was 8 seconds for $\mathcal{J}$ and also 8 seconds for $\mathcal{F}$. Thus, our saving of 124 seconds for $\mathcal{J}$ was made up of 132 saved in normal substring searches avoided less the extra 8 seconds on flags. Likewise, for $\mathcal{F}$, there were 338 seconds of normal substring searches saved.

# 7 Change flags combined with signatures

There is no obvious relationship between pairs of relators which fail the compatibility test and between pairs which are unchanged. This means that we can expect to combine both tests to reduce the amount of explicit substring searching even further. We have done so and we obtained the following results.

For $\mathcal{J}$, 6314024 out of 6413797 pairs of relators previously tested were unchanged or incompatible. This left only 99773 pairs to test, a saving of 98% over normal testing, a reduction of 92% on signature methods, and a reduction of 78% on change flags alone. The time taken was 39 seconds, a saving of 75% on normal methods and 54% as against signatures. However it took an extra 22% over the time taken for change flags alone. The substring searches saved (in addition to those already saved by change flags alone) were not enough to compensate for the signature computations. For $\mathcal{F}$, 8419814 out of 8614464 pairs of relators previously tested were unchanged or incompatible. This left only 194650 pairs to test, a saving of 98% over naive testing, a reduction of 87% on signature methods, and a reduction of 75% on change flags alone. The time taken was 139 seconds, a saving of 69% on normal methods and 53% as against signatures. However it took an extra 29% over the time taken for change flags alone. Again, the extra substring searches saved did not compensate for the signature computations.

The additional savings from adding signature tests to change flags confirms the relative independence of the methods. However, almost as much signature computation is done for this absolute saving as is done for the much greater absolute saving when used without change flags. There is a small reduction in signature comparisons, but it is the signature computation which takes most of the time. For the absolute savings obtained in conjunction with change flags, the signature computation takes so long as to be detrimental to total time taken.

# 8 Conclusions

We have studied the substring searching component of presentation manipulation algorithms used in computational group theory. In particular we have investigated improved methods based on signatures, on change flags, and on a combination of both of these. Detailed experiments reveal that all of these provide significant performance improvements and that, of them, change flags perform best. Some future directions of research with other techniques are indicated.

# 9 Acknowledgements

The authors were partially supported by Australian Research Council grant A49030651. We are grateful to Volkmar Felsch for suggesting the use of time stamps on changed relators, which led on to change flags. Some of the experimental work described here was done by James Darwin.